\documentclass[12pt,a4paper]{amsart}\usepackage{amssymb,euscript,array,bbm}

\newtheorem{Theorem}[equation]{Theorem}\newcommand{\bthe}{\begin{Theorem}}\newcommand{\ethe}{\end{Theorem}}
\newtheorem{Lemma}[equation]{Lemma}\newcommand{\bel}{\begin{Lemma}}\newcommand{\eel}{\end{Lemma}}
\newtheorem{Proposition}[equation]{Proposition}\newcommand{\bprop}{\begin{Proposition}}\newcommand{\eprop}{\end{Proposition}}
\newtheorem{Corollary}[equation]{Corollary}\newcommand{\bcor}{\begin{Corollary}}\newcommand{\ecor}{\end{Corollary}}
\newtheorem{Definition}[equation]{Definition}\newcommand{\bed}{\begin{Definition}}\newcommand{\eed}{\end{Definition}}
\newtheorem{Definition-Proposition}[equation]{Definition-Proposition}

\def\bpr{~\\{\em Proof.\ }}\newcommand{\epr}{$\bull$\\}
\newtheorem{Remark}[equation]{Remark}\newcommand{\beR}{\begin{Remark}\rm}\newcommand{\eeR}{\end{Remark}}
\newtheorem{Example}[equation]{Example}\newcommand{\bex}{\begin{Example}\rm}\newcommand{\eex}{\end{Example}}

\def\RR{{R}}

\newcommand{\beq}{\begin{equation}}\newcommand{\eeq}{\end{equation}}
\newcommand{\bem}{\begin{displaymath}}\newcommand{\eem}{\end{displaymath}}
\newcommand{\beqa}{\begin{eqnarray}}\newcommand{\eeqa}{\end{eqnarray}}
\newcommand{\bee}{\begin{enumerate}}\newcommand{\eee}{\end{enumerate}}
\newcommand{\bei}{\begin{itemize}}\newcommand{\eei}{\end{itemize}}
\newcommand{\bet}{\begin{tabular}{cccccccc}}\newcommand{\eet}{\end{tabular}}
\newcommand{\bpm}{\begin{pmatrix}}\newcommand{\epm}{\end{pmatrix}}
\newcommand{\bM}{\begin{matrix}}\newcommand{\eM}{\end{matrix}}
\newcommand{\ber}{\begin{array}{l}}\newcommand{\eer}{\end{array}}

\def\bull{\vrule height .9ex width .9ex depth -.1ex }

\def\li{~\\ $\bullet$ }

\def\cA{{\mathcal{A}}}\def\cB{{\mathcal{B}}}

\def\cU{\mathcal{U}}\def\cV{\mathcal{V}}
\def\cm{{\frak m}}

\def\one{{1\hspace{-0.1cm}\rm I}}
\def\zero{\mathbb{O}}

\def\C{\mathbb{C}}
\def\k{\mathbbm{k}}

\def\R{\mathbb{R}}\def\R{\mathbb{R}}

\def\X{\mathbb{X}}\def\Z{\mathbb{Z}}

\def\de{\delta}

  \def\tU{\tilde{U}}\def\tV{{\tilde{V}}}

\def\suml{\sum\limits}\def\oplusl{\mathop\oplus\limits}

\def\prodl{\prod\limits}

\def\sset{\subset}\def\sseteq{\subseteq}

\def\Mat{Mat(m,n;\RR)}\def\Matm{Mat(m,n;\cm)}

\newcommand{\bece}{\begin{center}}\newcommand{\eece}{\end{center}}
\newcommand{\quotients}[2]{{\footnotesize\left.\raisebox{0.4ex}{$#1$}\! / \!\raisebox{-0.4ex}{$#2$}\right.}}

\title{B\MakeLowercase{lock-diagonalization of matrices over local rings} II.}
\author{D\MakeLowercase{mitry} K\MakeLowercase{erner}}
\address{Department of Mathematics, Ben Gurion University of the Negev, P.O.B. 653, Be'er Sheva 84105, Israel.}
\email{dmitry.kerner@gmail.com}
\thanks{I was supported by the grant FP7-People-MCA-CIG, 334347.
 }
\date{\today}

\begin{document}
\maketitle  \setcounter{secnumdepth}{4}\setcounter{tocdepth}{1}
\begin{abstract}
Consider rectangular matrices over a local ring $R$. In the previous work we have obtained criteria for block-diagonalization of such matrices, i.e. $UAV=A_1\oplus A_2$, where $U,V$ are invertible matrices over $R$.
In this short note we extend the criteria to the decomposability of quiver representations over $R$.
\end{abstract}
%
%
%
%
%

\section{Introduction}
This work is the continuation of \cite{Kerner-Vinnikov.Part1}.
\subsection{Setup}
Let $(\RR,\cm)$ be a local (commutative, associative) ring over a field $\k$ of characteristic 0.
 As the simplest examples one can consider {\em regular} rings, e.g. formal power series, $\k[[x_1,\dots,x_p]]$, rational functions that are regular
at the origin, $\k[x_1,..,x_p]_{(\cm)}$, convergent power series, $k\{x_1,..,x_p\}$, when $\k$ is a normed field.
(If $\k=\R$ or $\k=\C$, one can consider the rings of germs of continuous or smooth functions as well.)
Usually we assume the ring to be non-Artinian, i.e. of positive Krull dimension (though $\RR$ can be not pure dimensional).

Denote by $\Mat$ the set of matrices with entries in $\RR$. In this paper we always assume: $1<m\le n$.
 Usually we assume that the matrices  "vanish at the origin", $A|_0=\zero$, i.e. $A\in\Matm$. Various matrix equivalences are important.
\li In commutative algebra matrices are considered up to the left-right equivalence, $A\stackrel{G_{lr}}{\sim}UAV$, where $(U,V)\in G_{lr}:=GL(m,\RR)\times GL(n,\RR)$.
\li In representation theory (of algebras/groups) the matrices are considered up to the conjugation, $A\stackrel{G_{conj}}{\sim}UAU^{-1}$, $U\in GL(m,\RR)$.
\li For the study of bilinear/quadratic/skew-symmetric forms one considers the congruence, $A\stackrel{G_{congr}}{\sim}UAU^T$, $U\in GL(n,\RR)$. (Note that we consider the non-primitive forms, i.e. $A$ vanishes mod $\cm$.)
\li More generally, one studies the matrix problems/representations of quivers. Each such representation consists of a collection of (rectangular) matrices and a prescribed transformation equivalence.

\

Unlike the case of classical linear algebra (over a field), the matrices over a ring cannot be diagonalized or
 brought to some nice/simple/canonical form.  For a given group action $G\circlearrowright\Mat$ the natural weaker question is the decomposability:
\beq
\text{\em Which matrices are block-diagonalizable, i.e. $A\stackrel{G}{\sim}\bpm A_1&\zero\\\zero&A_2\epm$?}
\eeq

In \cite{Kerner-Vinnikov.Part1} we have addressed this question for $G_{lr}$-equivalence.
Recall that the Fitting ideals (the ideals of $j\times j$ minors, $\{I_j(A)\}$) are invariant under $G_{lr}$-equivalence, \cite[\S20]{Eisenbud-book}. Thus the ideal of maximal minors of a block-diagonalizable matrix necessarily factorizes: $I_m(A)=I_{m_1}(A_1)I_{m_2}(A_2)$. We have obtained the following necessary and sufficient conditions for block-diagonalizability.
\bthe\label{Thm.Decomposability.L-R.}
1. \cite[Theorem 2]{Kerner-Vinnikov.Part1} (the case of square matrices)
\\Let $A\in Mat(m,m;\cm)$, $m>1$, with $\det(A)=f_1f_2$. Suppose each $f_i\in\RR$ is neither invertible
 nor a zero divisor and $f_1,f_2$ are relatively prime, i.e. $(f_1)\cap(f_2)=(f_1f_2)$.
 Then $A\stackrel{G_{lr}}{\sim} A_1\oplus A_2$, with $\det(A_i)=f_i$ iff $I_{m-1}(A)\sseteq (f_1)+(f_2)\sset \RR$.
\\2. \cite[Theorem 4]{Kerner-Vinnikov.Part1} (the case of rectangular matrices)
\\Let $A\in\Matm$, $m\le n$.
Suppose the ideal $I_m(A)$ does not annihilate any non-zero element of $\RR$, i.e. $ann_\RR(I_m(A))=\{0\}$.
  Suppose further that $ker(A)\sseteq I_m(A)\RR^{\oplus n}$.
  Suppose $I_m(A)=J_1J_2$, where
the (nontrivial) ideals $J_1,J_2\sset \RR$ are mutually prime, i.e. $J_1\cap J_2=J_1J_2$.
Then $A\stackrel{G_{lr}}{\sim} A_1\oplus A_2$ with $I_{m_i}(A_i)=J_i$ iff $I_{m-1}(A)\sseteq (J_1+J_2)\sset\RR$.
\ethe
(In the second case, when speaking of $ker(A)$, we consider $A$ as the map of free modules $\RR^{\oplus n}\stackrel{A}{\to}\RR^{\oplus m}$. The condition $ker(A)\sseteq I_m(A)\RR^{\oplus n}$ is the genericity assumption.)

\

 In this short note we extend the decomposability criteria to other equivalences. More precisely, we reduce the decomposability of quiver representations to the decomposability of matrices under $G_{lr}$-equivalence. This  reduction comes at the expense of enlarging the ring. However theorem \ref{Thm.Decomposability.L-R.} is "insensitive" to the dimension of the ring.  Thus the theorem can be used effectively  to obtain explicit decomposability criteria for various quivers/matrix problems.

\section{Quiver representations over local rings}
Given a quiver $Q$ with the set of vertices $I$. A representation of $Q$ over a ring $\RR$ is the collection of $R$-modules, $\{M_i\}_{i\in I}$ and of their morphisms $M_j\stackrel{A_{ij}}{\to}M_i$.

The following elementary observation is frequently used.
Given a ring $\RR$, consider the formal extension by some new variables, $\RR[[\{x_i\}]]$. Given two matrices, $A,B$ over $\RR[[\{x_i\}]]$ whose entries are linear in $\{x_i\}$, i.e. $A=\suml_i x_i A_i$ for some matrices $A_i$ over $\RR$. If $A=\tU B \tV$, where $\tU$, $\tV$ are invertible, with entries in $\RR[[\{x_i\}]]$, then there exist  invertible matrices $U,V$, over $\RR$, satisfying $A=UBV$. This leads to:
\bprop\label{Remark.Merging.Arrows.in.Quiver}
Fix a subgroup $G\sseteq G_{lr}$. Two tuples of matrices over $\RR$ are simultaneously $G$-equivalent, $(A_1,\dots,A_N)\stackrel{G}{\sim}(B_1,\dots,B_N)$, iff the corresponding matrices $\suml_i A_i x_i$,  $\suml_i B_i x_i$ are $G$-equivalent (over $\RR[[\{x_i\}]]$).
\eprop
We use this property for the groups $G_{lr}$, $G_{conj}$, $G_{congr}$.

\subsection{Replacing the quiver by a complete reduced quiver}
Fix a quiver $Q$ and its representation $\{A_{ij}\}$.
\li We can (and will) always assume that for any two vertices $(i,j)$ of $Q$ the quiver has arrows in both directions, $i\leftrightarrows j$. (Add all the missing arrows to the initial quiver and assume that the corresponding morphisms are zeros.)
\li We can (and will) assume that there are precisely two arrows: $i\to j$ and $j\to i$. If there are more, e.g. there is a tuple of morphisms $M_i\stackrel{(A^{(1)}_{ij},\dots,A^{(N)}_{ij})}{\to}M_j$, then we extend the ring to $\RR[[\{x_i\}]]$, and replace this tuple by one morphism $\suml_k A^{(k)}_{ij}x_k$. By proposition \ref{Remark.Merging.Arrows.in.Quiver}  we get an equivalent problem.

We call the so obtained quiver "a complete reduced" quiver.

\subsection{Embedding the quiver representations into representation of the Kronecker quiver}
Given a complete reduced quiver $Q_{\RR}$, with its representation $\{A_{ij}\}$, over $\RR$. Associate to it the following square matrix:
\beq
\{A_{ij}\}_{\substack{1\le i,j\le m}}\rightsquigarrow \cA_Q:=\bpm x_{11}A_{11}+y_1\one&x_{12}A_{12}&\cdots&x_{1m}A_{1m}\\
x_{21}A_{21}&x_{22}A_{22}+y_2\one&x_{23}A_{23}&\cdots
\\\cdots&\cdots&\cdots&\cdots
\\x_{m1}A_{m1}&\cdots&\cdots&x_{mm}A_{mm}+y_m\one\epm
\eeq
Here $\{x_{ij}\}$ and $\{y_i\}$ are some new formal variables (so that $\cA_Q$ is linear in these variables), while $\one$ denote the identity matrices of the appropriate sizes.
We consider $\cA_Q$ as a representation of the Kronecker quiver over $\RR[[\{x_{ij}\},\{y_i\}]]$.
\bprop\label{Thm.Embedding.Quiver.Rep}
Two quiver representations $\{A_{ij}\}$, $\{B_{ij}\}$ are equivalent (over $\RR$) iff the corresponding matrices $\cA_Q$ and $\cB_Q$ are $G_{lr}$-equivalent (over $\RR[[\{x_{ij}\},\{y_i\}]]$).
\eprop
Apparently this embedding is well known but we know the reference only for the particular case of conjugation, when $Q$ has just one vertex and one arrow, e.g. \cite[\S8]{Rao}.
\\\bpr
$\Rrightarrow$ The equivalence of representations means $A_{ij}=U_i B_{ij}U^{-1}_j$, for some matrices $\{U_i\in GL(m_i,\RR)\}$. Therefore
$\bpm U_1&\zero&\cdots\\\cdots&&\cdots\\\zero&\cdots&U_m\epm\cA_Q \bpm U^{-1}_1&\zero&\cdots\\\cdots&&\cdots\\\zero&\cdots&U^{-1}_m\epm=\cB_Q$.

$\Lleftarrow$ Suppose $\cB=\tilde\cU\cA_Q\tilde\cV$, where $\tilde\cU,\tilde\cV\in GL(\suml_i m_i,\RR[[\{x_{ij}\},\{y_i\}]])$. Recall that $\cB_Q$, $\cA_Q$ are linear in $\{x_{ij}\},\{y_i\}$. Thus, by proposition \ref{Remark.Merging.Arrows.in.Quiver}, we can choose invertible matrices $\cU$, $\cV$ over $\RR$ satisfying $\cB_Q=\cU\cA_Q\cV$. In the later equality consider the $y_i$ parts for each $i$. In particular, put $y_1=\cdots=y_m$ to get $\cU\cV=\one$.
   Now consider each $y_i$ separately to get: $\cU=\oplusl_i U_i$, $\cV=\oplusl_i U^{-1}_i$. Therefore
   $(\oplusl_i U_i)\cA_Q(\oplusl_i U^{-1}_i)=\cB_Q$. This implies $U_i A_{ij}U^{-1}_{ij}=B_{ij}$.
\epr

\subsection{Decomposability of $\{A_{ij}\}$ vs decomposability of $\cA_Q$}
\bed
$\cA_Q$ is quiver-block-diagonalizable if $\cA_Q\stackrel{G_{lr}}{\sim}\cA_1\oplus\cA_2$, where each of the determinants $det(\cA_k)$ contains a monomial $\prodl^m_{i=1}y^{l_i}_i$, where $0<l_i<m_i$.
\eed
\bprop\label{Thm.Reduction.quiver.rep.Decompos.to.A.decompos}
The representation $\{A_{ij}\}$ is decomposable iff the matrix $\cA_Q$ is quiver-block-diagonalizable.
\eprop
\bpr
$\Rrightarrow$ Suppose $\{A_{ij}\}$ is (non-trivially) decomposable, i.e. $\{U_i A_{ij}U^{-1}_j=\bpm A^{(1)}_{ij}&\zero\\\zero&A^{(2)}_{ij}\epm\}$. Denote $\cU=\oplus U_i$. Then $\cU\cA_Q\cU^{-1}$ consists of blocks, the $ij$'th block being $\bpm A^{(1)}_{ij}+\de_{ij}y_i\one&\zero\\\zero&A^{(2)}_{ij}+\de_{ij}y_i\one\epm$. (Here $\de_{ij}=1$ if $i=j$ and $0$ if $i\neq j$.) Note that each matrix $y_i\one$ is of non-zero size.

Then, by row/column permutations one can bring $\cU\cA_Q\cU^{-1}$ to the form $\cA_1\oplus\cA_2$, where each $\cA_k$ consists of blocks $A^{(k)}_{ij}+\de_{ij}y_i\one$. As each block $y_i\one$ is of non-zero size, $det(\cA_k)$ contains the monomial $\prodl^m_{i=1}y^{l_i}_i$ with $0<l_i$.

\

$\Lleftarrow$ Consider the images of $\cA_k$ under the projection $R\stackrel{\phi}{\to}\quotients{R}{\{x_{ij}\}}$.
 Note that $\phi(\cA_Q)$ is diagonal and its entries are linear in $\{y_i\}$, $\phi(\cA)=\oplus y_i C_i$ and $\sum rank(C_i)=rank(\cA_Q)$. Thus we can assume that the entries of each diagonal matrix $\phi(\cA_k)$ are linear in $\{y_i\}$ and the similar rank decomposition holds: $\phi(\cA_k)=\oplus y_i C^{(k)}_i$. Now $\phi(\cA_Q)$ and $\phi(\cA_1)\oplus \phi(\cA_2)$ are related by row/column permutations. Fix the corresponding matrices $U,V$, over $\k$, such that
  $U\phi(\cA_Q)V^{-1}=\bpm \phi(\cA_1)&\zero\\\zero&\phi(\cA_2)\epm$.

Using the embedding $\k\hookrightarrow\RR$ consider $U,V$ as matrices over $\RR$. Consider the matrix $U^{-1}\bpm\cA_1&\zero\\\zero&\cA_2\epm V$. It is of the form $\bpm y_1\one&\zero\\\dots&\dots&\dots\\\zero&& y_m\one\epm+\X$, where the block structure of $\X$ is:
 $\X_{ij}=\bpm \X^{(1)}_{ij}&\zero\\\zero&\X^{(2)}_{ij}\epm$ with $size(\X^{(k)}_{ij})=size(C^{(k)}_i)$.
Thus the matrix $U^{-1}\bpm\cA_1&\zero\\\zero&\cA_2\epm V$ comes from a decomposable representaton of the initial quiver. Now, invoke proposition \ref{Thm.Embedding.Quiver.Rep} to get: $\{A_{ij}\}$ is decomposable.
\epr
Combining this result with theorem \ref{Thm.Decomposability.L-R.} we get:
\bcor
1. If the representation $\{A_{ij}\}$ is decomposable then $det(\cA)\in\RR[[\{x_{ij}\},\{y_i\}]]$ is reducible.
\\2. Suppose $det(\cA)=f_1f_2$, where each of $f_1,f_2$ is a non-zero divisor and contains a monomial $\prodl_i y^{l_i}_i$. Suppose they are relatively prime, i.e. $(f_1)\cap(f_2)=(f_1f_2)$. Then the representation $\{A_{ij}\}$ is decomposable iff $I_{m-1}(\cA)\sset(f_1)+(f_2)$.
\ecor
In many examples $\det(\cA)$ is square free, thus $f_1,f_2$ are necessarily relatively prime. So, the corollary gives a very simple and effective decomposability criterion.

\section{Examples}
\subsection{Conjugation}
The conjugation, corresponds to the quiver with just one vertex (and several arrows). Explicitly, we are given a tuple of square matrices up to the simultaneous conjugation, $(A_1,\dots,A_m)\to U(A_1,\dots,A_m)U^{-1}$.
 Introduce the new variables, $x_1,\dots,x_m$, $y$. The associated matrix is $\cA_Q=\suml_i x_i A_i+y\one$. So that
  $\{A_i\}\sim\{B_i\}$ iff $\cA_Q\sim\cB_Q$. Proposition \ref{Thm.Reduction.quiver.rep.Decompos.to.A.decompos} reads:
\beq\ber
\text{the representation $\{A_i\}$ is $G_{conj}$-decomposable iff $\cA_Q=\suml_i x_iA_i+y\one$ is $G_{lr}$-block-diagonalizable}.
\eer\eeq
 As the simplest case consider the $2\times2$ matrices over $\RR$.
\bcor
Given $A=\bpm a_{11}&a_{12}\\a_{21}&a_{22}\epm\in Mat(2,2;\RR)$, suppose $tr^2(A)\neq 4det(A)$. Then $A$ is $G_{conj}$-decomposable iff $tr^2(A)-4det(A)$ is a full square in $\RR$ and moreover
 the elements $a_{12}$, $a_{21}$, $(a_{11}-a_{22})$ all belong to the ideal $\Big(\sqrt{tr^2(A)-4det(A)}\Big)$.
\ecor
\bpr
We want to check the $G_{lr}$-decomposability of $xA+y\one$.
  First of all we get that $det(xA+y\one)=y^2+xy\cdot tr(A)+x^2\cdot det(A)$ must factor over $\RR[[y]]$.
 Consider this as a quadratic equation for $y$. The roots are $y_{\pm}=x\frac{-tr(A)\pm\sqrt{tr^2(A)-4det(A)}}{2}$.
 Thus $tr^2(A)-4det(A)$ must be a full square in $\RR$.

To use theorem \ref{Thm.Decomposability.L-R.} we want the factors $(y-y_-)$, $(y-y_+)$ of
 $det(xA+y\one)$ to be relatively prime. As $y,x$ are independent variables, the factors are relatively prime iff $y_-\neq y_+$, i.e. $tr^2(A)-4det(A)\neq0$.

 Finally, if $tr^2(A)-4det(A)\neq0$, but is a full square in $\RR$, then the condition $I_1(xA+y\one)\sseteq J_1+J_2$ reads: $(xa_{12},xa_{21},2y+xtr(A),x(a_{11}-a_{22}))\sseteq (2y+xtr(A),x\sqrt{tr^2(A)-4det(A)})$. Which means $(a_{12},a_{21},a_{11}-a_{22})\sseteq (\sqrt{tr^2(A)-4det(A)})$.
\epr
 \bex
Let $A=\bpm a_{11}x_2&a_{12}x^k_1\\a_{21}x^l_1&a_{22}x_2\epm$, where $0\neq a_{ij}\in\k$ and $\k=\bar\k$ and $x_1,x_2$ are algebraically independent.
Then $tr^2(A)-4det(A)=(a_{11}-a_{22})^2x^2_2+a_{12}a_{21}x^{k+l}_1$.
If this expression is a full square then either $a_{11}=a_{22}$ and $k+l\in2\Z$ or $a_{12}a_{21}=0$.
\li Suppose $a_{11}=a_{22}$ and $k+l\in2\Z$. To ensure $tr^2(A)-4det(A)\neq0$ we assume $a_{11}a_{22}\neq0$. Then the condition  $\{a_{12}x^k_1,a_{21}x^l_1\}\sset (\sqrt{a_{12}a_{21}x^{k+l}_1})$ means: $k=l$.
\li Suppose $a_{12}a_{21}=0$. Then the condition  $\{a_{12}x^k_1,a_{21}x^l_1,(a_{11}-a_{22})x_2\}\sset ((a_{11}-a_{22})x_2)$ means: $a_{12}=a_{21}=0$, i.e. $A$ is already diagonal.
\li Finally, if $a_{11}=a_{22}=0$, but $a_{12}a_{21}\neq0$ then the matrix is not $G_{conj}$-diagonalizable, e.g. by checking the characteristic polynomial.

Summarizing: if $A$ is not diagonal then it is $G_{conj}$-diagonalizable iff $k=l$,  $a_{11}=a_{22}=1$.
\eex
\bex
 Similarly consider $A=\bpm x_2&x^k_1&0\\0&x_2&x^l_1\\-x^{3n-k-l}_1&0&x_2\epm$. Then $det(A)=(x_2-x^n_1)(x^2_2+x_2x^n_1+x^{2n}_1)$, note that the two factors are mutually prime. Thus the $G_{lr}$-decomposability holds iff $k=l=n$.
  To check the $G_{conj}$-decomposability we consider $\cA=\bpm y+x_2&x^n_1&0\\0&y+x_2&x^n_1\\-x^n_1&0&y+x_2\epm$. Then $I_2(\cA)=(y+x_2)^2+((y+x_2)x^n_1)+(x^{2n}_1)\sseteq (x_2-x^n_1)+(x^2_2+x_2x^n_1+x^{2n}_1)$. Thus $A$ is $G_{conj}$-decomposable.
\eex

\subsection{The "star" quiver} Consider the one-vertex quiver, with $l$ incoming arrows and $m$ outgoing arrows,
$Q=\bM A_{1+1}\nwarrow\uparrow\nearrow A_k\\\hspace{1cm}\bullet\circlearrowleft C\\B_1\nearrow\uparrow\nwarrow B_l\eM$. The corresponding matrix is
\beq
\cA_Q=\bpm x_{00}C+y_0\one&\zero&\dots&\zero&x_{0,l+1}A_{l+1}&\dots&x_{0,k}A_k\\x_{10}B_1&y_1\one&\zero\\
\dots&&&\dots\\x_{l0}B_l&\zero&\dots&y_l\one&\zero&\dots\\\zero&\dots&\dots&\zero&y_{l+1}\one&\zero
\\\dots&&&&&\dots&y_k\one\epm
\eeq
Here $det(\cA_Q)=det(x_{00}C+y_0\one)\prodl^k_{i=1}y^{m_i}_i$, so any factorization $det(\cA_Q)=f_1f_2$  is impossible
with $f_1,f_2$ relatively prime. Theorem \ref{Thm.Decomposability.L-R.} does not produce any decomposability criterion here.

\subsubsection{The string quiver}
Consider the quiver
$Q=\bullet\ber A_1\\\leftrightarrows\\B_1\eer\bullet\dots\bullet\dots\bullet\ber A_{m-1}\\\leftrightarrows\\B_{m-1}\eer\bullet$. The corresponding matrix is:
\beq
\cA_Q=\bpm y_1\one&x_{12}A_1&\zero&\dots&\zero\\x_{21}B_1&y_2\one&x_{23}A_2&\zero
\\\zero&x_{32}B_2&y_3\one&x_{34}A_3&\zero\\\dots&\dots&\dots&\dots&x_{m-1,m}A_{m-1}\\
 \zero&\dots&\zero&x_{m,m-1}B_{m-1}&y_m\one\epm
\eeq
The decomposability of a representation of $Q$ is controlled by $\cA_Q$.
\bex Consider the simplest case, $Q=\bullet\ber A\\\leftrightarrows\\B\eer\bullet$, where $A,B\in Mat(2,2;\RR)$.  Then $\det\cA_Q=y^2_1y^2_2-y_1y_2x_{12}x_{21}tr(AB)+x^2_{12}x^2_{21}det(AB)$. Thus the decomposability implies:
$tr^2(AB)-4det(AB)$ is a full square in $\RR$. Further, if $\cA_Q$ is decomposable then
the generators of $I_3(\cA_Q)$, i.e. the elements
\beq
y_1y^2_2,\ y^2_1y_2,\ \{y_1y_2x_{12}A_{ij}\},\ \{y_1y_2x_{21}B_{ij}\},\ x^2_{12}x_{21}det(A)\{B_{ij}\},\
x_{12}x^2_{21}det(B)\{A_{ij}\}
\eeq
 belong to the ideal: $(2y_1y_2-x_{12}x_{21}tr(AB))+\Big(x_{12}x_{21}\sqrt{tr^2(AB)-4det(AB)}\Big)$. As $y_1,y_2$, $x_{12},x_{21}$ are algebraically independent, only two cases are possible:
\bei
\item $tr(AB)=0$ and the elements $det(A)\{B_{ij}\}$, $det(B)\{A_{ij}\}$ belong to the ideal $(\sqrt{-det(AB)})$.
In this case the factors of $\det\cA_Q$ are mutually prime iff $det(AB)\neq0$.
 If $det(AB)\neq0$ then theorem \ref{Thm.Decomposability.L-R.} ensures decomposability.
If $det(AB)=0$, then the condition on $det(A)\{B_{ij}\}$, $det(B)\{A_{ij}\}$ forces: either $A=\zero$ or $B=\zero$ or $det(A)=0=det(B)$.
\item $det(AB)=a\cdot tr^2(AB)$, where $a-\frac{1}{4}$ is invertible in $\RR$. Then the factors of $\det\cA_Q$ are mutually prime and together generate the ideal $(y_1y_2)+(x_{12}x_{21}tr(AB))$.
In this case the conditions give: both $det(A)\{B_{ij}\}$  $det(B)\{A_{ij}\}$ belong to the ideal $(tr(AB))$.
\eei
\eex

\end{document}